\documentclass[12pt,oneside]{article}
\usepackage{amsmath,amssymb,amsfonts,amsthm}
\usepackage{color}
\usepackage[all]{xy}
\usepackage{graphicx}
\usepackage{color}
\usepackage{makeidx}
\usepackage{psfrag}
\usepackage{epstopdf}
\usepackage{multirow}
\usepackage{rotating}
\usepackage{maplestd2e}
\newcommand {\N}{\mathcal{N}}

\def\Section#1{\vspace{30truept}\addtocounter{section}{1}\setcounter{thm}{0}\setcounter{equation}{0}
{\noindent\Large\bf\arabic{section}.#1}\par \vspace{12pt}}

\newtheorem{thm}{Theorem}[section]

\newtheorem{defn}[thm]{Definition}

\newtheorem{rem}[thm]{Remark}

\numberwithin{equation}{section}
\newcommand\undersym[2]{\raisebox{-7pt}{\tiny$#2$}{\kern-8pt}\mbox{$#1$}}
\newcommand\undersymm[2]{\raisebox{-8pt}{\tiny$#2$}{\kern-15pt}\mbox{$#1$}}
\newcommand\overast[1]{\raisebox{9pt}{\small$\ast$}{\kern-9pt}\mbox{$#1$}}
\newcommand\overlind[1]{\raisebox{10pt}{\small$\overline{{\hspace{2pt}}\star}$}{\kern-7.5pt}\mbox{$#1$}}
\newcommand\overlinc[1]{\raisebox{10pt}{\tiny$\overline{{\hspace{2pt}}\circ}$}{\kern-7.5pt}\mbox{$#1$}}
\newcommand\overlina[1]{\raisebox{10pt}{\small$\overline{{\hspace{1pt}}\ast}$}{\kern-7.5pt}\mbox{$#1$}}
\newcommand\overcirc[1]{\raisebox{10pt}{\tiny{$\circ$}}{\kern-7.5pt}\mbox{$#1$}}
\newcommand\overdiamond[1]{\raisebox{10pt}{\small$\star$}{\kern-7.5pt}\mbox{$#1$}}

\newcommand\tovercirc[1]{\raisebox{5pt}{\tiny{$\circ$}}{\kern-5.5pt}\mbox{$#1$}}
\newcommand\toverdiamond[1]{\raisebox{5pt}{\tiny$\star$}{\kern-5.5pt}\mbox{$#1$}}
\newcommand\toverast[1]{\raisebox{5pt}{\tiny$\ast$}{\kern-5pt}\mbox{$#1$}}

\begin{document}
\title{\bf{Computing nullity and kernel  vectors using NF-package: Counterexamples}}
%\footnote{ArXiv Number: }}
\author{\bf{ Nabil L. Youssef$^{\,1}$ and S. G.
Elgendi$^{2}$}}
\date{}
%\thanks{\it Department of Mathematics, etc}
%\pagestyle{fancy}
             % End of preamble and beginning of text.
\maketitle                     % Produces the title.
\vspace{-1.16cm}
\begin{center}
{$^{1}$Department of Mathematics, Faculty of Science,\\ Cairo
University, Giza, Egypt}\\
\end{center}

\begin{center}
{$^{2}$Department of Mathematics, Faculty of Science,\\ Benha
University, Benha,
 Egypt}
\end{center}

\begin{center}
E-mails: nlyoussef@sci.cu.edu.eg, nlyoussef2003@yahoo.fr\\
{\hspace{1.8cm}}salah.ali@fsci.bu.edu.eg, salahelgendi@yahoo.com
\end{center}
\smallskip
\vspace{1cm} \maketitle
\smallskip
{\vspace{-1.1cm}} \noindent{\bf \begin{center}
Abstract\end{center}}
A computational technique for calculating  nullity vectors and kernel vectors, using the   new Finsler  package,   is introduced. As an application, three interesting counterexamples are given. The first counterexample shows   that the two distributions   $\mathrm{Ker}_R$ and  $\N_R$ do not coincide. The second  shows that  the nullity distribution $\N_{P^\circ}$  is not completely integrable. The third shows that the nullity distribution $\N_\mathfrak{R}$ is not a sub-distribution of  the nullity distribution  $\N_{R^\circ}$.

\vspace{7pt}

\newpage
%%%%%%%%%%%%%%%%%%%%%%%%%%%%%%%%%%%%%%%%%%%%%%%%%%%%%%%% %%%%%%%%%%%%%%%%%%%%%%%%%%%%%%%%%%%%%%%%%%%%%%%%%%%%%%%%%%%%%%%%%%%%%%%%%%%%%%%%%%%
%\vspace{1cm}
%Introduction
%\vspace{30truept}\centerline{\Large\bf{Introduction}}\vspace{12pt}
\begin{center}
\large \textbf{Introduction}
\end{center}

\par

In the applicable examples  of    Finsler geometry in mathematics, physics and the other branches of science, the calculations are often very tedious to perform. This  takes a lot of effort and time. So, we have to find an  alternative method  to do these calculations.   One of the benefits of using computer is the manipulation  of the  complicated calculations.   This  enables  to study various  examples in different dimensions in various applications  (cf., for example, \cite{r101}, \cite{Rutz2},  \cite{shen2005}, \cite{r93}, \cite{gamal},\cite{Portugal1},\cite{wanas}). The FINSLER package \cite{Rutz3} included in \cite{hbfinsler1} and the new Finsler package \cite{CFG} are   good  illustrations of using computer in the applications of Finsler geometry. \\

 In this paper, we  use the new Finsler (NF-) package \cite{CFG} to introduce a computational technique to calculate the components of  nullity vectors and  kernel vectors.  As an application of this method, we construct  three interesting counterexamples. The first shows that the kernel distribution  $\mathrm{Ker}_R$ and the  nullity distribution $\N_R$ associated with the h-curvature $R$ of Cartan connection do not coincide, in accordance with  \cite{ND-Zadeh}. The second proves that the nullity distribution $\N_{P^\circ}$ associated with the hv-curvature \, $\overcirc{P}$ of  Berwald connection is not completely integrable. Finally, the third counterexample  shows that the nullity distribution $\N_\mathfrak{R}$ associated with the curvature $\mathfrak{R}$ of  Barthel connection is not a sub-distribution of the nullity distribution $\N_{R^\circ}$ associated with the h-curvature \, $\overcirc{R}$ of  Berwald connection.\\

 Following the Klein-Grifone approach to Finsler geometry (\cite{r21}, \cite{r22}, \cite{r27}), let $(M,F)$ be a Finsler space, where $F$ is a Finsler structure defined on an $n$-dimensional smooth manifold $M$. Let H(TM) (resp. V(TM)) be the horizontal (resp. vertical) sub-bundle of the bundle TTM. We use the notations $R$  and $P$ for the h-curvature and  hv-curvature of Cartan connection respectively. We also use the notations\, $\overcirc{R}$ and\, $\overcirc{P}$ for the h-curvature and hv-curvature of Berwald connection respectively. Finally, $\mathfrak{R} $ will denote the curvature of the Cartan non-linear connection (Barthel connection).

%%%%%%%%%%%%%%%%%%%%%%%%%%%%%%%%%%%%%%%%%%%%%%%%%%%%%%%%%%%%%%%%%%%%%%%%%%%%%%%%%%%%%%%%%%%%%%%%%%%%%%%%%%%%%%%%%%%%%%%%%%%%%%%%%%%%%%%%%%%%%%%%%%%%%%%%%%%
\Section{ Nullity  and kernel vectors by  the NF-package }

In this section, we  use the  New Finsler (NF-) package \cite{CFG}, which is an extended and modified version of \cite{Rutz3}, to introduce a computational method for the calculation of  nullity vectors and  kernel vectors.

\begin{defn}\label{nr} Let $R$ be the h-curvature tensor of Cartan connection.
The nullity space of $R$ at a point $z\in TM$ is the subspace of $H_z(TM)$ defined by
$$\mathcal{N}_R(z):=\{X\in H_z(TM) : \,  R(X,Y)Z=0, \, \,\forall\, Y,Z\in H_z(TM)\}.$$
The dimension of $\mathcal{N}_R(z)$, denoted by $\mu_R(z)$, is the index of nullity of $R$ at $z$.

 If $\mu_R(z)$ is constant,  the map $\mathcal{N}_R:z\mapsto \mathcal{N}_R(z) $ defines a distribution $\mathcal{N}_R$ of
rank $\mu_R$ called nullity distribution of $R$.

Any  vector field belonging to the nullity distribution is called  a nullity vector field.
 \end{defn}

\begin{defn}\label{ker}
The kernel space  $\mathrm{Ker}_{R}(z)$  of the h-curvature  ${R}$ at a point $z\in TM$ is the subspace of $H_z(TM)$ defined by
$$\mathrm{Ker}_{R}(z)=\{X\in H_z(TM):  \, {R}(Y,Z)X=0, \, \forall\, Y,Z\in H_z(TM)\}. $$

As in Definition \ref{nr}, the map $z\mapsto \mathrm{Ker}_{R}(z) $ defines a distribution called the kernel   distribution of $R$.
Any  vector field belonging to the kernel distribution is called  a kernel vector field.
\end{defn}

To  calculate the nullity vectors and kernel vectors using the NF-package, let  us  recall  some
 instructions  to make  the use of this package easier.  When we write, for example,  N[i,-j] we  mean $N^i_j$, i.e.,   positive (resp. negative) index means that it is  contravariant  (resp.  covariant). To  lower or raise  an index by  the metric or the inverse  metric, just change its sign from  positive to negative or  vice versa.   The command  \lq\lq \emph{tdiff}(N[i,-j], X[k])\rq\rq means ${\partial}_k N^i_j$, the command  \lq\lq \emph{tddiff}(N[i,-j], Y[k])\rq\rq means $\dot{\partial}_k N^i_j$ and the command \lq\lq\emph{Hdiff}(N[i,-j], X[k])\rq\rq means $\delta_k N^i_j$. To introduce the definition of a tensor, we use the command \lq\lq\emph{definetensor}\rq\rq and to display its components, we use the command \lq\lq\emph{show}\rq\rq as will be seen soon.\\

Now, let $Z\in\N_R$ be a nullity vector. Then, $Z$ can be written locally  in the form $Z=Z^ih_i$, where $Z^i$  are the  components of the nullity vector $Z$ with respect to the basis   $\{h_i\}$  of the horizontal space, where  $h_i:=\frac{\partial}{\partial x^i}-N^j_i\frac{\partial}{\partial y^j}$ and $N^j_i$ are the coefficients of Barthel connection;  $i, j=1,...,n$.  The equation  ${R}(Z,X)Y=0$, $\forall\, X, Y\in H(TM)$, is written locally in the form
  $$Z^j{R}^i_{hjk}=0.$$
   To derive the resulting system from   $Z^j{R}^i_{hjk}=0$,   we first compute the components $R^i_{hjk}$ using the  NF-package. Then, we define  a new tensor by the command \lq\lq\emph{definetensor}\rq\rq as follows:
  \smallskip
  \begin{maplegroup}
\begin{mapleinput}
\mapleinline{active}{1d}{definetensor(RCZ[h,-i,-k] = RC[h,-i,-j,-k]*Z[j]); }{\[\]}
\end{mapleinput}
\end{maplegroup}
\begin{maplegroup}
\begin{mapleinput}
\mapleinline{active}{1d}{show(RCZ[h,-i,-k]); }{\[\]}
\end{mapleinput}
  \mapleresult
\end{maplegroup}
\smallskip
Putting  $RCZ[h,-i,-k]=0$, we obtain  a homogenous system of algebraic equations. Solving this system, we  get the components $Z^i$.
\begin{rem}
\em{It should be noted that we must  not use  the notation $X=X^ih_i$ nor the notation $Y=Y^ih_i$ for   nullity vectors because $RC[h,-i,-j,-k]*X[j]$ and $RC[h,-i,-j,-k]*Y[j]$  mean to Maple   $x^jR^h_{ijk}$ and $y^jR^h_{ijk}$ respectively,  which both are not the correct expressions for  nullity vectors.  }
\end{rem}

In a similar way, we compute the components of a  kernel vector. Let $W=W^ih_i\in \mathrm{Ker}_R$, then $R(X,Y)W=0$, $\forall\, X, Y\in H(TM)$. This locally  gives the homogenous system of algebraic equations:
$$W^h{R}^i_{hjk}=0.$$
Then by the  NF-package, we can define
\smallskip
\begin{maplegroup}
\begin{mapleinput}
\mapleinline{active}{1d}{definetensor(RCW[h,-j,-k] = RC[h,-i,-j,-k]*W[i]); }{\[\]}
\end{mapleinput}
\end{maplegroup}
\begin{maplegroup}
\begin{mapleinput}
\mapleinline{active}{1d}{show(RCW[h,-j,-k]); }{\[\]}
\end{mapleinput}
  \mapleresult
\end{maplegroup}
\smallskip
Putting $RCW[h,-j,-k]=0$  and  solving the resulting  system, we get the components $W^i$ of the kernel vector $W$.

 %%%%%%%%%%%%%%%%%%%%%%%%%%%%%%%%%%%%%%%%%%%%%%%%%%%%%%%%%%%%%%%%%%%%%%%%%%%%%%%%%%%%%%%%%%%%%%%%%%%%%%%%%%%%%%%%%%%%%%%%%%%%%%%%%%%%%%%%%%%%
\newpage
\Section{ Applications and counterexamples }

\par
In this section,  we provide  three interesting  counterexamples. We perform the computations using the above mentioned technique and the      NF-package. We also make use of the technique of simplification of tesor expressions \cite{CFG}.

\vspace{5pt}

The nullity distributions associated with Cartan connection are  studied in \cite{ND-cartan}. The following example shows that \emph{the nullity space $\N_R$ of the h-curvature $R$ of Cartan connection and  the kernel  $\mathrm{Ker}_R$ do not coincide.}

\bigskip

\noindent \textbf{Example 1}

 Let $M=\{(x^1,...,x^4)\in \mathbb{R}^4|\,  x^2>0\}$, $U=\{(x^1,...,x^4;y^1,...,y^4)\in \mathbb{R}^4 \times \mathbb{R}^4: \, y^2\neq~0~, y^4\neq 0\}\subset TM$.  Let  $F$ be defined on $U$ by
 $$F := \, ({{{{\it x2}}^{2}{{\it y1}}^{4}+{{\it y2}}^{4}+{{\it y3}}^{4}+{{\it y4}}^{4}}})^{1/4}.$$

By Maple program and NF-package we can perform the following calculations.

\bigskip

\begin{maplegroup}
\begin{mapleinput}
\mapleinline{active}{1d}{F0 := sqrt(x2\symbol{94}2*y1\symbol{94}4+y2\symbol{94}4+y3\symbol{94}4+y4\symbol{94}4);
}{}
\end{mapleinput}
\mapleresult
\begin{maplelatex}
\mapleinline{inert}{2d}{F0 := sqrt(x2^2*y1^4+y2^4+y3^4+y4^4)}{\[ F0\, := \, \sqrt{{{\it x2}}^{2}{{\it y1}}^{4}+{{\it y2}}^{4}+{{\it y3}}^{4}+{{\it y4}}^{4}\\
}\]}
\end{maplelatex}
\end{maplegroup}

\smallskip

\begin{maplegroup}
\begin{Maple Normal}{
\textbf{Barthel connection}}\end{Maple Normal}

\end{maplegroup}
\begin{maplegroup}
\begin{mapleinput}
\mapleinline{active}{1d}{show(N[i,-j]);
}{}
\end{mapleinput}
\mapleresult

\end{maplegroup}
\begin{maplelatex}
\mapleinline{inert}{2d}{}{}{\[
   N^{{\it x1} }_{{\it x1} }=\frac{1}{3}\frac{{\it y2} }{{\it x2} } \hspace{1cm} N^{{\it x1} }_{{\it x2} }=\frac{1}{3}\frac{{\it y1} }{{\it x2} } \hspace{1cm}
  N^{{\it x2} }_{{\it x1} }=-\frac{1}{3}\frac{{\it x2}{\it y1}^{3}}{{\it y2}^{2}} \hspace{1cm} N^{{\it x2} }_{{\it x2} }=\frac{1}{6}\frac{{\it x2}{\it y1}^{4}}{{\it y2}^{3}}\]}
\end{maplelatex}
\mapleresult

\vspace{7pt}

\begin{maplegroup}
\begin{Maple Normal}{
\textbf{ h-curvature R of Cartan connection }}\end{Maple Normal}

\end{maplegroup}
\begin{maplegroup}
\begin{mapleinput}
\mapleinline{active}{1d}{show(RC[h, -i, -j, -k]); }{\[\]}
\end{mapleinput}
\mapleresult
\begin{maplelatex}
\mapleinline{inert}{2d}{}{}{\[ {\it RC}^{{\it x1} }_{{\it x2 x1 x2} }=-\frac{1}{18}\frac{3{\it x2}^{4}{\it y1}^{8}+2{\it x2}^{2}{\it y1}^{4}{\it y4}^{4}+2{\it y3}^{4}{\it x2}^{2}{\it y1}^{4}+13{\it x2}^{2}{\it y1}^{4}{\it y2}^{4}+4{\it y2}^{8}+8{\it y3}^{4}{\it y2}^{4}+8{\it y2}^{4}{\it y4}^{4}}{{\it x2}^{2}\left({\it x2}^{2}{\it y1}^{4}+{\it y2}^{4}+{\it y3}^{4}+{\it y4}^{4}\right){\it y2}^{4}}\]}
\end{maplelatex}
\begin{maplelatex}
\mapleinline{inert}{2d}{}{}{\[ {\it RC}^{{\it x2} }_{{\it x1 x1 x2} }=\frac{1}{18}\frac{\left({\it x2}^{4}{\it y1}^{8}+2{\it y3}^{4}{\it x2}^{2}{\it y1}^{4}+2{\it x2}^{2}{\it y1}^{4}{\it y4}^{4}+7{\it x2}^{2}{\it y1}^{4}{\it y2}^{4}+8{\it y2}^{4}{\it y4}^{4}+8{\it y3}^{4}{\it y2}^{4}+12{\it y2}^{8}\right){\it y1}^{2}}{{\it y2}^{6}\left({\it x2}^{2}{\it y1}^{4}+{\it y2}^{4}+{\it y3}^{4}+{\it y4}^{4}\right)}\]}
\end{maplelatex}
\mapleresult
\begin{maplelatex}
\mapleinline{inert}{2d}{}{}{\[ {\it RC}^{{\it x1} }_{{\it x1 x1 x2} }=\frac{1}{9}\frac{{\it y1}^{3}\left(4{\it y2}^{4}+{\it x2}^{2}{\it y1}^{4}\right)}{\left({\it x2}^{2}{\it y1}^{4}+{\it y2}^{4}+{\it y3}^{4}+{\it y4}^{4}\right){\it y2}^{3}}\hspace{1cm} {\it RC}^{{\it x1} }_{{\it x3 x1 x2} }=\frac{1}{18}\frac{\left(4{\it y2}^{4}+{\it x2}^{2}{\it y1}^{4}\right){\it y3}^{3}}{{\it x2}^{2}{\it y2}^{3}\left({\it x2}^{2}{\it y1}^{4}+{\it y2}^{4}+{\it y3}^{4}+{\it y4}^{4}\right)}\]}
\end{maplelatex}
\mapleresult
\begin{maplelatex}
\mapleinline{inert}{2d}{}{}{\[ {\it RC}^{{\it x1} }_{{\it x4 x1 x2} }=\frac{1}{18}\frac{\left(4{\it y2}^{4}+{\it x2}^{2}{\it y1}^{4}\right){\it y4}^{3}}{{\it x2}^{2}{\it y2}^{3}\left({\it x2}^{2}{\it y1}^{4}+{\it y2}^{4}+{\it y3}^{4}+{\it y4}^{4}\right)}\hspace{.5cm} {\it RC}^{{\it x2} }_{{\it x2 x1 x2} }=-\frac{1}{9}\frac{{\it y1}^{3}\left(4{\it y2}^{4}+{\it x2}^{2}{\it y1}^{4}\right)}{\left({\it x2}^{2}{\it y1}^{4}+{\it y2}^{4}+{\it y3}^{4}+{\it y4}^{4}\right){\it y2}^{3}}\]}
\end{maplelatex}
\mapleresult
\begin{maplelatex}
\mapleinline{inert}{2d}{}{}{\[ {\it RC}^{{\it x2} }_{{\it x3 x1 x2} }=-\frac{1}{18}\frac{{\it y1}^{3}{\it y3}^{3}\left(4{\it y2}^{4}+{\it x2}^{2}{\it y1}^{4}\right)}{{\it y2}^{6}\left({\it x2}^{2}{\it y1}^{4}+{\it y2}^{4}+{\it y3}^{4}+{\it y4}^{4}\right)}\hspace{.6cm} {\it RC}^{{\it x2} }_{{\it x4 x1 x2} }=-\frac{1}{18}\frac{{\it y1}^{3}{\it y4}^{3}\left(4{\it y2}^{4}+{\it x2}^{2}{\it y1}^{4}\right)}{{\it y2}^{6}\left({\it x2}^{2}{\it y1}^{4}+{\it y2}^{4}+{\it y3}^{4}+{\it y4}^{4}\right)}\]}
\end{maplelatex}
\mapleresult
\begin{maplelatex}
\mapleinline{inert}{2d}{}{}{\[ {\it RC}^{{\it x3} }_{{\it x1 x1 x2}}=\frac{1}{18}\frac{\left(4{\it y2}^{4}+{\it x2}^{2}{\it y1}^{4}\right){\it y1}^{2}{\it y3} }{\left({\it x2}^{2}{\it y1}^{4}+{\it y2}^{4}+{\it y3}^{4}+{\it y4}^{4}\right){\it y2}^{3}}\hspace{1cm} {\it RC}^{{\it x3} }_{{\it x2 x1 x2} }=-\frac{1}{18}\frac{\left(4{\it y2}^{4}+{\it x2}^{2}{\it y1}^{4}\right){\it y3}{\it y1}^{3}}{\left({\it x2}^{2}{\it y1}^{4}+{\it y2}^{4}+{\it y3}^{4}+{\it y4}^{4}\right){\it y2}^{4}}\]}
\end{maplelatex}
\mapleresult
\begin{maplelatex}
\mapleinline{inert}{2d}{}{}{\[ {\it RC}^{{\it x4} }_{{\it x1 x1 x2} }=\frac{1}{18}\frac{\left(4{\it y2}^{4}+{\it x2}^{2}{\it y1}^{4}\right){\it y4}{\it y1}^{2}}{\left({\it x2}^{2}{\it y1}^{4}+{\it y2}^{4}+{\it y3}^{4}+{\it y4}^{4}\right){\it y2}^{3}}\hspace{1cm} {\it RC}^{{\it x4} }_{{\it x2 x1 x2} }=-\frac{1}{18}\frac{\left(4{\it y2}^{4}+{\it x2}^{2}{\it y1}^{4}\right){\it y4}{\it y1}^{3}}{\left({\it x2}^{2}{\it y1}^{4}+{\it y2}^{4}+{\it y3}^{4}+{\it y4}^{4}\right){\it y2}^{4}}\]}
\end{maplelatex}
\end{maplegroup}

\vspace{11pt}

\begin{maplegroup}
\begin{Maple Normal}{
\textbf{$R$-Nullity vectors}}\end{Maple Normal}

\end{maplegroup}

\begin{maplegroup}
\begin{mapleinput}
\mapleinline{active}{1d}{definetensor(RCW[h, -i, -k] = RC[h, -i, -j, -k]*W[j]); }{\[\]}
\end{mapleinput}
\end{maplegroup}

\begin{mapleinput}
\mapleinline{active}{1d}{show(RCW[h, -i, -k]); }{\[\]}
\end{mapleinput}
\mapleresult
\begin{maplelatex}
\mapleinline{inert}{2d}{}{}{\[ {\it RCW}^{{\it x1} }_{{\it x2 x1} }=\frac{1}{18}\frac{\left(3{\it x2}^{4}{\it y1}^{8}+13{\it x2}^{2}{\it y1}^{4}{\it y2}^{4}+2{\it x2}^{2}{\it y1}^{4}{\it y4}^{4}+2{\it y3}^{4}{\it x2}^{2}{\it y1}^{4}+8{\it y2}^{4}{\it y4}^{4}+4{\it y2}^{8}+8{\it y3}^{4}{\it y2}^{4}\right)W^{{\it x2} }}{{\it x2}^{2}\left({\it x2}^{2}{\it y1}^{4}+{\it y2}^{4}+{\it y3}^{4}+{\it y4}^{4}\right){\it y2}^{4}}\]}
\end{maplelatex}
\mapleresult
\begin{maplelatex}
\mapleinline{inert}{2d}{}{}{\[ {\it RCW}^{{\it x1} }_{{\it x2 x2} }=-\frac{1}{18}\frac{\left(3{\it x2}^{4}{\it y1}^{8}+13{\it x2}^{2}{\it y1}^{4}{\it y2}^{4}+2{\it x2}^{2}{\it y1}^{4}{\it y4}^{4}+2{\it y3}^{4}{\it x2}^{2}{\it y1}^{4}+8{\it y2}^{4}{\it y4}^{4}+4{\it y2}^{8}+8{\it y3}^{4}{\it y2}^{4}\right)W^{{\it x1} }}{{\it x2}^{2}\left({\it x2}^{2}{\it y1}^{4}+{\it y2}^{4}+{\it y3}^{4}+{\it y4}^{4}\right){\it y2}^{4}}\]}
\end{maplelatex}
\mapleresult
\begin{maplelatex}
\mapleinline{inert}{2d}{}{}{\[ {\it RCW}^{{\it x2} }_{{\it x1 x1} }=-\frac{1}{18}\frac{\left({\it x2}^{4}{\it y1}^{8}+7{\it x2}^{2}{\it y1}^{4}{\it y2}^{4}+2{\it y3}^{4}{\it x2}^{2}{\it y1}^{4}+2{\it x2}^{2}{\it y1}^{4}{\it y4}^{4}+12{\it y2}^{8}+8{\it y2}^{4}{\it y4}^{4}+8{\it y3}^{4}{\it y2}^{4}\right){\it y1}^{2}W^{{\it x2} }}{{\it y2}^{6}\left({\it x2}^{2}{\it y1}^{4}+{\it y2}^{4}+{\it y3}^{4}+{\it y4}^{4}\right)}\]}
\end{maplelatex}
\mapleresult
\begin{maplelatex}
\mapleinline{inert}{2d}{}{}{\[ {\it RCW}^{{\it x2} }_{{\it x1 x2} }=\frac{1}{18}\frac{\left({\it x2}^{4}{\it y1}^{8}+7{\it x2}^{2}{\it y1}^{4}{\it y2}^{4}+2{\it y3}^{4}{\it x2}^{2}{\it y1}^{4}+2{\it x2}^{2}{\it y1}^{4}{\it y4}^{4}+12{\it y2}^{8}+8{\it y2}^{4}{\it y4}^{4}+8{\it y3}^{4}{\it y2}^{4}\right){\it y1}^{2}W^{{\it x1} }}{{\it y2}^{6}\left({\it x2}^{2}{\it y1}^{4}+{\it y2}^{4}+{\it y3}^{4}+{\it y4}^{4}\right)}\]}
\end{maplelatex}
\mapleresult
\begin{maplelatex}
\mapleinline{inert}{2d}{}{}{\[ {\it RCW}^{{\it x1} }_{{\it x1  x1} }=-\frac{1}{9}\frac{{\it y1}^{3}\left(4{\it y2}^{4}+{\it x2}^{2}{\it y1}^{4}\right)W^{{\it x2} }}{\left({\it x2}^{2}{\it y1}^{4}+{\it y2}^{4}+{\it y3}^{4}+{\it y4}^{4}\right){\it y2}^{3}}\hspace{1cm} {\it RCW}^{{\it x1} }_{{\it x1  x2} }=\frac{1}{9}\frac{{\it y1}^{3}\left(4{\it y2}^{4}+{\it x2}^{2}{\it y1}^{4}\right)W^{{\it x1} }}{\left({\it x2}^{2}{\it y1}^{4}+{\it y2}^{4}+{\it y3}^{4}+{\it y4}^{4}\right){\it y2}^{3}}\]}
\end{maplelatex}
\mapleresult
\begin{maplelatex}
\mapleinline{inert}{2d}{}{}{\[ {\it RCW}^{{\it x1} }_{{\it x3 x1} }=-\frac{1}{18}\frac{\left(4{\it y2}^{4}+{\it x2}^{2}{\it y1}^{4}\right){\it y3}^{3}W^{{\it x2} }}{{\it x2}^{2}{\it y2}^{3}\left({\it x2}^{2}{\it y1}^{4}+{\it y2}^{4}+{\it y3}^{4}+{\it y4}^{4}\right)} \hspace{.3cm}
{\it RCW}^{{\it x1} }_{{\it x3x2} }=\frac{1}{18}\frac{\left(4{\it y2}^{4}+{\it x2}^{2}{\it y1}^{4}\right){\it y3}^{3}W^{{\it x1} }}{{\it x2}^{2}{\it y2}^{3}\left({\it x2}^{2}{\it y1}^{4}+{\it y2}^{4}+{\it y3}^{4}+{\it y4}^{4}\right)}\]}
\end{maplelatex}
\mapleresult
\begin{maplelatex}
\mapleinline{inert}{2d}{}{}{\[ {\it RCW}^{{\it x1} }_{{\it x4 x1} }=-\frac{1}{18}\frac{\left(4{\it y2}^{4}+{\it x2}^{2}{\it y1}^{4}\right){\it y4}^{3}W^{{\it x2} }}{{\it x2}^{2}{\it y2}^{3}\left({\it x2}^{2}{\it y1}^{4}+{\it y2}^{4}+{\it y3}^{4}+{\it y4}^{4}\right)}\hspace{.3cm}
 {\it RCW}^{{\it x1} }_{{\it x4 x2} }=\frac{1}{18}\frac{\left(4{\it y2}^{4}+{\it x2}^{2}{\it y1}^{4}\right){\it y4}^{3}W^{{\it x1} }}{{\it x2}^{2}{\it y2}^{3}\left({\it x2}^{2}{\it y1}^{4}+{\it y2}^{4}+{\it y3}^{4}+{\it y4}^{4}\right)}\]}
\end{maplelatex}
\mapleresult
\begin{maplelatex}
\mapleinline{inert}{2d}{}{}{\[ {\it RCW}^{{\it x2} }_{{\it x2 x1} }=\frac{1}{9}\frac{{\it y1}^{3}\left(4{\it y2}^{4}+{\it x2}^{2}{\it y1}^{4}\right)W^{{\it x2} }}{\left({\it x2}^{2}{\it y1}^{4}+{\it y2}^{4}+{\it y3}^{4}+{\it y4}^{4}\right){\it y2}^{3}}\hspace{1.2cm}
 {\it RCW}^{{\it x2} }_{{\it x2 x2} }=-\frac{1}{9}\frac{{\it y1}^{3}\left(4{\it y2}^{4}+{\it x2}^{2}{\it y1}^{4}\right)W^{{\it x1} }}{\left({\it x2}^{2}{\it y1}^{4}+{\it y2}^{4}+{\it y3}^{4}+{\it y4}^{4}\right){\it y2}^{3}}\]}
\end{maplelatex}
\mapleresult
\begin{maplelatex}
\mapleinline{inert}{2d}{}{}{\[ {\it RCW}^{{\it x2} }_{{\it x3 x1} }=\frac{1}{18}\frac{{\it y1}^{3}{\it y3}^{3}\left(4{\it y2}^{4}+{\it x2}^{2}{\it y1}^{4}\right)W^{{\it x2} }}{{\it y2}^{6}\left({\it x2}^{2}{\it y1}^{4}+{\it y2}^{4}+{\it y3}^{4}+{\it y4}^{4}\right)}\hspace{1.1cm}
{\it RCW}^{{\it x2} }_{{\it x3 x2} }=-\frac{1}{18}\frac{{\it y1}^{3}{\it y3}^{3}\left(4{\it y2}^{4}+{\it x2}^{2}{\it y1}^{4}\right)W^{{\it x1} }}{{\it y2}^{6}\left({\it x2}^{2}{\it y1}^{4}+{\it y2}^{4}+{\it y3}^{4}+{\it y4}^{4}\right)}\]}
\end{maplelatex}
\mapleresult
\begin{maplelatex}
\mapleinline{inert}{2d}{}{}{\[ {\it RCW}^{{\it x2} }_{{\it x4 x1} }=\frac{1}{18}\frac{{\it y1}^{3}{\it y4}^{3}\left(4{\it y2}^{4}+{\it x2}^{2}{\it y1}^{4}\right)W ^{{\it x2} }}{{\it y2}^{6}\left({\it x2}^{2}{\it y1}^{4}+{\it y2}^{4}+{\it y3}^{4}+{\it y4}^{4}\right)}\hspace{1cm}
 {\it RCW}^{{\it x2} }_{{\it x4 x2} }=-\frac{1}{18}\frac{{\it y1}^{3}{\it y4}^{3}\left(4{\it y2}^{4}+{\it x2}^{2}{\it y1}^{4}\right)W^{{\it x1} }}{{\it y2}^{6}\left({\it x2}^{2}{\it y1}^{4}+{\it y2}^{4}+{\it y3}^{4}+{\it y4}^{4}\right)}\]}
\end{maplelatex}
\mapleresult
\begin{maplelatex}
\mapleinline{inert}{2d}{}{}{\[ {\it RCW}^{{\it x3} }_{{\it x1 x1} }=-\frac{1}{18}\frac{\left(4{\it y2}^{4}+{\it x2}^{2}{\it y1}^{4}\right){\it y3} {\it y1}^{2}W^{{\it x2} }}{\left({\it x2}^{2}{\it y1}^{4}+{\it y2}^{4}+{\it y3}^{4}+{\it y4}^{4}\right){\it y2}^{3}}\hspace{1cm}
 {\it RCW}^{{\it x3} }_{{\it x1  x2} }=\frac{1}{18}\frac{\left(4{\it y2}^{4}+{\it x2}^{2}{\it y1}^{4}\right){\it y3} {\it y1}^{2}W^{{\it x1} }}{\left({\it x2}^{2}{\it y1}^{4}+{\it y2}^{4}+{\it y3}^{4}+{\it y4}^{4}\right){\it y2}^{3}}\]}
\end{maplelatex}
\mapleresult
\begin{maplelatex}
\mapleinline{inert}{2d}{}{}{\[ {\it RCW}^{{\it x3} }_{{\it x2 x1} }=\frac{1}{18}\frac{\left(4{\it y2}^{4}+{\it x2}^{2}{\it y1}^{4}\right){\it y3} {\it y1}^{3}W^{{\it x2} }}{\left({\it x2}^{2}{\it y1}^{4}+{\it y2}^{4}+{\it y3}^{4}+{\it y4}^{4}\right){\it y2}^{4}}\hspace{1cm}
 {\it RCW}^{{\it x3} }_{{\it x2  x2} }=-\frac{1}{18}\frac{\left(4{\it y2}^{4}+{\it x2}^{2}{\it y1}^{4}\right){\it y3} {\it y1}^{3}W^{{\it x1} }}{\left({\it x2}^{2}{\it y1}^{4}+{\it y2}^{4}+{\it y3}^{4}+{\it y4}^{4}\right){\it y2}^{4}}\]}
\end{maplelatex}
\mapleresult
\begin{maplelatex}
\mapleinline{inert}{2d}{}{}{\[ {\it RCW}^{{\it x4} }_{{\it x1  x1} }=-\frac{1}{18}\frac{\left(4{\it y2}^{4}+{\it x2}^{2}{\it y1}^{4}\right){\it y1}^{2}{\it y4} W^{{\it x2} }}{\left({\it x2}^{2}{\it y1}^{4}+{\it y2}^{4}+{\it y3}^{4}+{\it y4}^{4}\right){\it y2}^{3}}\hspace{1cm}
 {\it RCW}^{{\it x4} }_{{\it x1  x2} }=\frac{1}{18}\frac{\left(4{\it y2}^{4}+{\it x2}^{2}{\it y1}^{4}\right){\it y1}^{2}{\it y4} W^{{\it x1} }}{\left({\it x2}^{2}{\it y1}^{4}+{\it y2}^{4}+{\it y3}^{4}+{\it y4}^{4}\right){\it y2}^{3}}\]}
\end{maplelatex}
\mapleresult
\begin{maplelatex}
\mapleinline{inert}{2d}{}{}{\[ {\it RCW}^{{\it x4} }_{{\it x2  x1} }=\frac{1}{18}\frac{\left(4{\it y2}^{4}+{\it x2}^{2}{\it y1}^{4}\right){\it y4} {\it y1}^{3}W^{{\it x2} }}{\left({\it x2}^{2}{\it y1}^{4}+{\it y2}^{4}+{\it y3}^{4}+{\it y4}^{4}\right){\it y2}^{4}}\hspace{1cm}
 {\it RCW}^{{\it x4} }_{{\it x2  x2} }=-\frac{1}{18}\frac{\left(4{\it y2}^{4}+{\it x2}^{2}{\it y1}^{4}\right){\it y4} {\it y1}^{3}W^{{\it x1} }}{\left({\it x2}^{2}{\it y1}^{4}+{\it y2}^{4}+{\it y3}^{4}+{\it y4}^{4}\right){\it y2}^{4}}\]}
\end{maplelatex}
\mapleresult
\vspace{7pt}
Putting ${\it RCW}^{\it h }_{\it ij }=0$, then we have a system of algebraic equations. The NF-package yields the following solution: $W^1=W^2=0,   W^3=s, W^4=t,; \,s,t\in \mathbb{R}$. Then,  any  nullity vector $W$ has the form
\begin{equation}\label{null.1}
W=sh_3+th_4.
\end{equation}

%\vspace{2pt}

\begin{maplegroup}
\begin{mapleinput}
\mapleinline{active}{2d}{}{\[\]}
\end{mapleinput}
\end{maplegroup}
\begin{Maple Normal}{
\begin{Maple Normal}{
\textbf{$R$-Kernel vectors}}\end{Maple Normal}

}\end{Maple Normal}

\begin{maplegroup}
\begin{mapleinput}
\mapleinline{active}{1d}{definetensor(RCZ[h, -j, -k] = RC[h, -i, -j, -k]*Z[i]); }{\[\]}
\end{mapleinput}
\end{maplegroup}
\begin{maplegroup}
\begin{mapleinput}
\mapleinline{active}{1d}{show(RCZ[h, -j, -k]); }{\[\]}
\end{mapleinput}
\mapleresult
\begin{maplelatex}
\mapleinline{inert}{2d}{}{\[ {\it RCZ}^{{\it x1} }_{{\it x1  x2} }=\frac{1}{9}\frac{{\it y1}^{3}\left(4{\it y2}^{4}+{\it x2}^{2}{\it y1}^{4}\right)Z^{{\it x1} }}{\left({\it x2}^{2}{\it y1}^{4}+{\it y2}^{4}+{\it y3}^{4}+{\it y4}^{4}\right){\it y2}^{3}}\]}
\mapleinline{inert}{2d}{}{\[\hspace{1.7cm}-\frac{1}{18}\frac{\left(3{\it x2}^{4}{\it y1}^{8}+13{\it x2}^{2}{\it y1}^{4}{\it y2}^{4}+2{\it x2}^{2}{\it y1}^{4}{\it y4}^{4}+2{\it y3}^{4}{\it x2}^{2}{\it y1}^{4}+8{\it y2}^{4}{\it y4}^{4}+4{\it y2}^{8}+8{\it y3}^{4}{\it y2}^{4}\right)Z^{{\it x2} }}{{\it x2}^{2}\left({\it x2}^{2}{\it y1}^{4}+{\it y2}^{4}+{\it y3}^{4}+{\it y4}^{4}\right){\it y2}^{4}}\]}
\mapleinline{inert}{2d}{}{\[\hspace{1.7cm}+\frac{1}{18}\frac{\left(4{\it y2}^{4}+{\it x2}^{2}{\it y1}^{4}\right){\it y3}^{3}Z^{{\it x3} }}{{\it x2}^{2}{\it y2}^{3}\left({\it x2}^{2}{\it y1}^{4}+{\it y2}^{4}+{\it y3}^{4}+{\it y4}^{4}\right)}+\frac{1}{18}\frac{\left(4{\it y2}^{4}+{\it x2}^{2}{\it y1}^{4}\right){\it y4}^{3}Z^{{\it x4} }}{{\it x2}^{2}{\it y2}^{3}\left({\it x2}^{2}{\it y1}^{4}+{\it y2}^{4}+{\it y3}^{4}+{\it y4}^{4}\right)}\]}
\end{maplelatex}
\mapleresult
\begin{maplelatex}
\mapleinline{inert}{2d}{}{\[ {\it RCZ}^{{\it x2} }_{{\it x1  x2} }=\frac{1}{18}\frac{\left({\it x2}^{4}{\it y1}^{8}+7{\it x2}^{2}{\it y1}^{4}{\it y2}^{4}+2{\it y3}^{4}{\it x2}^{2}{\it y1}^{4}+2{\it x2}^{2}{\it y1}^{4}{\it y4}^{4}+12{\it y2}^{8}+8{\it y2}^{4}{\it y4}^{4}+8{\it y3}^{4}{\it y2}^{4}\right){\it y1}^{2}Z^{{\it x1} }}{{\it y2}^{6}\left({\it x2}^{2}{\it y1}^{4}+{\it y2}^{4}+{\it y3}^{4}+{\it y4}^{4}\right)}\]}
\mapleinline{inert}{2d}{}{\[\hspace{1.5cm}-\frac{1}{9}\frac{{\it y1}^{3}\left(4{\it y2}^{4}+{\it x2}^{2}{\it y1}^{4}\right)Z^{{\it x2} }}{\left({\it x2}^{2}{\it y1}^{4}+{\it y2}^{4}+{\it y3}^{4}+{\it y4}^{4}\right){\it y2}^{3}}-\frac{1}{18}\frac{{\it y1}^{3}{\it y3}^{3}\left(4{\it y2}^{4}+{\it x2}^{2}{\it y1}^{4}\right)Z^{{\it x3} }}{{\it y2}^{6}\left({\it x2}^{2}{\it y1}^{4}+{\it y2}^{4}+{\it y3}^{4}+{\it y4}^{4}\right)}-\frac{1}{18}\frac{{\it y1}^{3}{\it y4}^{3}\left(4{\it y2}^{4}+{\it x2}^{2}{\it y1}^{4}\right)Z^{{\it x4} }}{{\it y2}^{6}\left({\it x2}^{2}{\it y1}^{4}+{\it y2}^{4}+{\it y3}^{4}+{\it y4}^{4}\right)}\]}
\end{maplelatex}
\mapleresult
\begin{maplelatex}
\mapleinline{inert}{2d}{}{\[ {\it RCZ}^{{\it x3} }_{{\it x1  x2} }=\frac{1}{18}\frac{\left(4{\it y2}^{4}+{\it x2}^{2}{\it y1}^{4}\right){\it y3} {\it y1}^{2}Z^{{\it x1} }}{\left({\it x2}^{2}{\it y1}^{4}+{\it y2}^{4}+{\it y3}^{4}+{\it y4}^{4}\right){\it y2}^{3}}-\frac{1}{18}\frac{\left(4{\it y2}^{4}+{\it x2}^{2}{\it y1}^{4}\right){\it y3} {\it y1}^{3}Z^{{\it x2} }}{\left({\it x2}^{2}{\it y1}^{4}+{\it y2}^{4}+{\it y3}^{4}+{\it y4}^{4}\right){\it y2}^{4}}\]}
\end{maplelatex}
\mapleresult
\begin{maplelatex}
\mapleinline{inert}{2d}{}{\[ {\it RCZ}^{{\it x4} }_{{\it x1  x2} }=\frac{1}{18}\frac{\left(4{\it y2}^{4}+{\it x2}^{2}{\it y1}^{4}\right){\it y1}^{2}{\it y4} Z^{{\it x1} }}{\left({\it x2}^{2}{\it y1}^{4}+{\it y2}^{4}+{\it y3}^{4}+{\it y4}^{4}\right){\it y2}^{3}}-\frac{1}{18}\frac{\left(4{\it y2}^{4}+{\it x2}^{2}{\it y1}^{4}\right){\it y4} {\it y1}^{3}Z^{{\it x2} }}{\left({\it x2}^{2}{\it y1}^{4}+{\it y2}^{4}+{\it y3}^{4}+{\it y4}^{4}\right){\it y2}^{4}}\]}
\end{maplelatex}
\mapleresult
\end{maplegroup}

\bigskip

Putting ${\it RCZ}^{\it h }_{\it ij }=0$, we obtain a system of algebraic equations. The NF-package yields the solution:
 $Z^1=\frac{sy_1}{y_2}$, $Z^2=s$, $Z^3=t$ and  $Z^4=\frac{s(x_2y_1^4+y_2^4+2y_3^4+2y_4^4)-ty_2y_3^3}{y_2y_4^3}$.

 Then, any  kernel vector $Z$ should have  the form
 \begin{equation}\label{ker1}Z=s\left(\frac{y_1}{y_2}h_1+h_2+\frac{x_2y_1^4+y_2^4+2y_3^4+2y_4^4}{y_2y_4^3}h_4\right)+t\left(h_3-\frac{y_3^3}{y_4^3}h_4\right).
 \end{equation}
(for simplicity, we have written $x_i$ and $y_i$ instead of $x^i$ and $y^i$ respectively)

 Comparing (\ref{null.1}) and (\ref{ker1}), we find   no values for $s$ and $t$  which  make $Z=W$. Consequently,
$\N_R$ and $\mathrm{Ker}_R$ can not coincide.

\bigskip

In \cite{Nabil.1} Youssef proved that the nullity distribution $\N_{R^\circ}$ associated with the h-curvature \, $\overcirc{R}$ of Berwald connection is completely integrable. He conjectured that the nullity distribution $\N_{P^\circ}$  of  the hv-curvature \,$\overcirc{P}$ of Berwald connection is not completely integrable. In the next example, \emph{we show that his conjecture is true}.

\bigskip

\noindent \textbf{Example 2}

 Let $M=\mathbb{R}^3$, $U=\{(x^1,x^2,x^3;y^1,y^2,y^3)\in \mathbb{R}^3 \times \mathbb{R}^3: \,y^1\neq 0\}\subset TM$.  Let  $F$ be defined on $U$ by
 $$F := \,{{\rm e}^{-{\it x1}}} \left( {{\it y2}}^{3}+{{\rm e}^{-{\it x1x3}}}{\it y3}\,{{\it y1}}^{2} \right)^{1/3}.$$

By Maple program and NF-package, we can perform the following calculations.

\bigskip

\begin{maplegroup}
\begin{mapleinput}
\mapleinline{active}{1d}{F0 := exp(-2*x1)*(y2\symbol{94}3+exp(-x1*x3)*y3*y1\symbol{94}2)\symbol{94}(2/3);
}{}
\end{mapleinput}
\mapleresult
\begin{maplelatex}
\mapleinline{inert}{2d}{F0 := exp(-2*x1)*(y2^3+exp(-x1*x3)*y3*y1^2)^(2/3)}{\[ {\it F0}\, := \,{{\rm e}^{-{\it 2 x1}}} \left( {{\it y2}}^{3}+{{\rm e}^{-{\it x1x3}}}{\it y3}\,{{\it y1}}^{2} \right)^{2/3}\]}
\end{maplelatex}
\end{maplegroup}

\vspace{7pt}

\begin{maplegroup}
\begin{Maple Normal}{
\textbf{Barthel connection}}\end{Maple Normal}

\end{maplegroup}
\begin{maplegroup}
\begin{mapleinput}
\mapleinline{active}{1d}{show(N[i,-j]);
}{}
\end{mapleinput}
\mapleresult
\begin{maplelatex}
\mapleinline{inert}{2d}{}{}{\[
   N^{{\it x1} }_{{\it x1} }=-\frac{1}{2}\left(3+{\it x3} \right){\it y1} \hspace{1cm} N^{{\it x2} }_{{\it x1} }=-\frac{3}{4}{\it y2}\hspace{1cm}
   N^{{\it x2} }_{{\it x2} }=-\frac{3}{4}{\it y1}\]}
\end{maplelatex}
\mapleresult
   \begin{maplelatex}
\mapleinline{inert}{2d}{}{}{\[
    N^{{\it x3} }_{{\it x1} }=-\frac{3}{4}\frac{{\it y2}^{3}}{{\it y1}^{2}e^{-x1 x3 }}\hspace{1cm}
   N^{{\it x3} }_{{\it x2} }=\frac{9}{4}\frac{{\it y2}^{2}}{{\it y1} e^{-x1 x3 }} \hspace{2cm} N^{{\it x3} }_{{\it x3} }=-{\it y3} {\it x1}\]}
\end{maplelatex}
\mapleresult
\end{maplegroup}
\vspace{7pt}

\begin{maplegroup}
\begin{Maple Normal}{
\textbf{ hv-curvature  \,$\overcirc{P}$ of Berwald connection}}\end{Maple Normal}
\end{maplegroup}
\begin{maplegroup}
\begin{mapleinput}
\mapleinline{active}{1d}{show(PB[h,-i,-j,-k]);  }{\[\]}
\end{mapleinput}
\begin{maplelatex}
\mapleinline{inert}{2d}{}{}{\[
  {\it PB}^{{\it x3} }_{{\it x1 x1 x1} }=-\frac{9}{2}\frac{{\it y2}^{3}}{{\it y1}^{4}e^{-x1 x3 }} \hspace{1cm}  {\it PB}^{{\it x3} }_{{\it x1 x1 x2} }=\frac{9}{2}\frac{{\it y2}^{2}}{{\it y1}^{3}e^{-x1 x3 }} \]}
\end{maplelatex}
\mapleresult
  \begin{maplelatex}
\mapleinline{inert}{2d}{}{}{\[
  {\it PB}^{{\it x3} }_{{\it x1 x2 x2} }=-\frac{9}{2}\frac{{\it y2} }{{\it y1}^{2}e^{-x1 x3 }} \hspace{1cm} {\it PB}^{{\it x3} }_{{\it x2 x2 x2} }=\frac{9}{2{\it y1} e^{-x1 x3 }}\]}
\end{maplelatex}
\mapleresult
\end{maplegroup}

\vspace{7pt}
\begin{maplegroup}
\begin{Maple Normal}{
\textbf{ $\overcirc{P}$-Nullity vectors}}\end{Maple Normal}

\end{maplegroup}

\begin{maplegroup}
\begin{mapleinput}
\mapleinline{active}{1d}{definetensor(PBW[i,-h,-k] = PB[i,-h,-j,-k]*W[j]); }{\[\]}
\end{mapleinput}
\end{maplegroup}
\begin{maplegroup}
\begin{mapleinput}
\mapleinline{active}{1d}{show(PBW[i,-h,-k]); }{\[\]}
\end{mapleinput}
\mapleresult
\begin{maplelatex}
\mapleinline{inert}{2d}{}{\[ {\it PBW}^{{\it x3} }_{{\it x1 x1} }=-\frac{9}{2}\frac{{\it y2}^{3}W^{{\it x1} }}{{\it y1}^{4}e^{-x1 x3 }}+\frac{9}{2}~\frac{{\it y2}^{2}W^{{\it x2} }}{{\it y1}^{3}e^{-x1 x3 }}\hspace{1cm}
{\it PBW}^{{\it x3} }_{{\it x2 x2} }=-\frac{9}{2}\frac{{\it y2}W^{{\it x1} }}{{\it y1}^{2}e^{-x1 x3 }}+\frac{9}{2}~\frac{{\it}W^{{\it x2} }}{{\it y1}e^{-x1 x3 }}\]}
\end{maplelatex}
\mapleresult
\begin{maplelatex}
\mapleinline{inert}{2d}{}{\[ {\it PBW}^{{\it x3} }_{{\it x1 x2} }=\frac{9}{2}\frac{{\it y2}^{2}W^{{\it x1} }}{{\it y1}^{3}e^{-x1 x3 }}-\frac{9}{2}~\frac{{\it y2}W^{{\it x2} }}{{\it y1}^{2}e^{-x1 x3 }}\hspace{1.5cm}
{\it PBW}^{{\it x3} }_{{\it x2 x1} }=\frac{9}{2}\frac{{\it y2}^{2}W^{{\it x1} }}{{\it y1}^{3}e^{-x1 x3 }}-\frac{9}{2}~\frac{{\it y2}W^{{\it x2} }}{{\it y1}^{2}e^{-x1 x3 }}\]}
\end{maplelatex}
\mapleresult
\end{maplegroup}

\bigskip

Putting ${\it PBW}^{{\it h} }_{{\it ij} }=0$, we get a system of algebraic equations. We have two cases:

 \noindent The first case  is $y2=0$ and  the solution in this case  is $W^1=s$, $W^2=0$ and $W^3=t$. Hence, any $\overcirc{P}$-nullity vector is written in the form  $W=sh_1+ th_3$. Take two nullity vectors  $X,Y \in \N_{P^\circ}$ such that $X=h_1$ and  $Y=h_3$. Their Lie bracket $[X,Y]=-\frac{y_1}{2}\frac{\partial}{\partial y_1}+y_3\frac{\partial}{\partial y_3}$,  which is vertical.

 \noindent The second case is $y2\neq 0$ and the solution in this case is  $W^1=s$, $W^2=\frac{y_2}{y_1}s$ and $W^3=t$. Then any  \, $\overcirc{P}$-nullity vector is written in the form  $W=s(h_1+\frac{y_2}{y_1}h_2)+ th_3$. Let $X$ and $Y$ be the  two nullity vectors in $\N_{P^\circ}$ given by $X=h_1+\frac{y_2}{y_1}h_2$ and $Y=h_3$. By computing their Lie  bracket, we find that  $[X,Y]=-\frac{y_1}{2}\frac{\partial}{\partial y_1}+y_3\frac{\partial}{\partial y_3}$,  which is vertical.

 \noindent  Consequently, in both cases the Lie   bracket $[X,Y]$ does not  belong to $\N_{P^\circ}$.\\

Let   $\N_{R^\circ}$  and $\N_\mathfrak{R}$ be the  nullity distributions associated with the h-curvature \, $\overcirc{R}$ of Berwald connection and  the curvature $\mathfrak{R}$ of the Barthel connection respectively. In \cite{Nabil.2}, Youssef proved that $\N_{R^\circ} \subseteq \N_\mathfrak{R}$. The following example  shows that\emph{ the converse is not true:  that is $\N_{R^\circ}$ is a proper sub-distribution of $\N_\mathfrak{R}$}.

\bigskip
\noindent \textbf{Example 3}

\smallskip

Let $M= \mathbb{R}^4$, $U=\{(x^1,\cdots,x^4;y^1,\cdots,y^4)\in \mathbb{R}^4 \times \mathbb{R}^4: \, y^2\neq 0, \,y^4\neq 0 \}\subset TM$.  Let  $F$ be defined on $U$ by
 $$F := \left(\,{{\rm e}^{-{\it x2}}}{\it y1}\,\sqrt [3]{{{\it y2}}^{3}+{{\it y3}}^{3}+{{\it y4}}^{3}}\right)^{1/2}.$$
By Maple program and NF-package, we can perform the following calculations.

\bigskip

\begin{maplegroup}
\begin{mapleinput}
\mapleinline{active}{1d}{F0 := exp(-x2)*y1*(y2\symbol{94}3+y3\symbol{94}3+y4\symbol{94}3)\symbol{94}(1/3);
}{}
\end{mapleinput}
\mapleresult
\begin{maplelatex}
\mapleinline{inert}{2d}{}{\[ {\it F0}\, :=  \,{{\rm e}^{-{\it x2}}}{\it y1}\,\sqrt [3]{{{\it y2}}^{3}+{{\it y3}}^{3}+{{\it y4}}^{3}}\]}
\end{maplelatex}
\end{maplegroup}

\vspace{7pt}

\begin{maplegroup}
\begin{Maple Normal}{
\textbf{Barthel connection}}\end{Maple Normal}

\end{maplegroup}
\begin{maplegroup}
\begin{mapleinput}
\mapleinline{active}{1d}{show(N[i,-j]);
}{}
\end{mapleinput}
\mapleresult
\end{maplegroup}
\begin{maplelatex}
\mapleinline{inert}{2d}{}{}{\[
  N^{{\it x2} }_{{\it x2} }=-\frac{1}{4}\frac{4{\it y2}^{3}+{\it y3}^{3}+{\it y4}^{3}}{{\it y2}^{2}}\hspace{1cm}
 N^{{\it x2} }_{{\it x3} }=\frac{3}{4}\frac{{\it y3}^{2}}{{\it y2} }  \hspace{2cm}
   N^{{\it x2} }_{{\it x4} }=\frac{3}{4}\frac{{\it y4}^{2}}{{\it y2} }\]}
\end{maplelatex}
 \mapleresult
   \begin{maplelatex}
\mapleinline{inert}{2d}{}{}{\[
 N^{{\it x3} }_{{\it x2} }=-\frac{3}{4}{\it y3} \hspace{1cm}
  N^{{\it x3} }_{{\it x3} }=-\frac{3}{4}{\it y2} \hspace{1cm}
 N^{{\it x4} }_{{\it x2} }=-\frac{3}{4}{\it y4}\hspace{1cm}
  N^{{\it x4} }_{{\it x4} }=-\frac{3}{4}{\it y2}\]}
\end{maplelatex}
 \mapleresult

\vspace{7pt}

\begin{maplegroup}
\begin{Maple Normal}{
\textbf{Curvature $\mathfrak{R}$ of the Barthel connection }}\end{Maple Normal}

\end{maplegroup}
\begin{maplegroup}
\begin{mapleinput}
\mapleinline{active}{1d}{show(RG[i, -j, -k]); }{\[\]}
\end{mapleinput}
\begin{maplelatex}
\mapleinline{inert}{2d}{}{}{\[
 {\it RG}^{{\it x2} }_{{\it x2 x3} }=-\frac{3}{16}\frac{{\it y3}^{2}\left({\it y2}^{3}+{\it y3}^{3}+{\it y4}^{3}\right)}{{\it y2}^{4}}\hspace{1cm}
 {\it RG}^{{\it x3} }_{{\it x2 x3} }=\frac{3}{16}\frac{{\it y2}^{3}+{\it y3}^{3}+{\it y4}^{3}}{{\it y2}^{2}}\]}
\end{maplelatex}
 \mapleresult
\begin{maplelatex}
\mapleinline{inert}{2d}{}{}{\[
{\it RG}^{{\it x2} }_{{\it x2 x4} }=-\frac{3}{16}\frac{{\it y4}^{2}\left({\it y2}^{3}+{\it y3}^{3}+{\it y4}^{3}\right)}{{\it y2}^{4}}\hspace{1cm}
 {\it RG}^{{\it x4} }_{{\it x2 x4} }=\frac{3}{16}\frac{{\it y2}^{3}+{\it y3}^{3}+{\it y4}^{3}}{{\it y2}^{2}}\]}
\end{maplelatex}
 \mapleresult
 \begin{maplelatex}
\mapleinline{inert}{2d}{}{}{\[
 {\it RG}^{{\it x3} }_{{\it x3 x4} }=\frac{9}{16}\frac{{\it y4}^{2}}{{\it y2} }\hspace{3.5cm}
 {\it RG}^{{\it x4} }_{{\it x3 x4} }=-\frac{9}{16}\frac{{\it y3}^{2}}{{\it y2} }\]}
\end{maplelatex}
 \mapleresult

\end{maplegroup}
\vspace{7pt}
\begin{maplegroup}
\begin{Maple Normal}{
\textbf{$\mathfrak{R}$-nullity vectors}}\end{Maple Normal}

\end{maplegroup}

\begin{maplegroup}
\begin{mapleinput}
\mapleinline{active}{1d}{definetensor(RGZ[i, -j] = RG[i, -j, -k]*Z[k]); }{\[\]}
\end{mapleinput}
\end{maplegroup}
\begin{maplegroup}
\begin{mapleinput}
\mapleinline{active}{1d}{show(RGZ[i, -j]); }{\[\]}
\end{mapleinput}
\mapleresult
 \begin{maplelatex}
\mapleinline{inert}{2d}{}{}{\[{\it RGZ}^{{\it x2} }_{{\it x2} }=-\frac{3}{16}\frac{{\it y3}^{2}\left({\it y2}^{3}+{\it y3}^{3}+{\it y4}^{3}\right)Z^{{\it x3} }}{{\it y2}^{4}}-\frac{3}{16}\frac{{\it y4}^{2}\left({\it y2}^{3}+{\it y3}^{3}+{\it y4}^{3}\right)Z^{{\it x4} }}{{\it y2}^{4}}\]}
\end{maplelatex}
 \mapleresult
 \begin{maplelatex}
\mapleinline{inert}{2d}{}{}{\[{\it RGZ}^{{\it x2} }_{{\it x3} }=\frac{3}{16}\frac{{\it y3}^{2}\left({\it y2}^{3}+{\it y3}^{3}+{\it y4}^{3}\right)Z^{{\it x2} }}{{\it y2}^{4}}\hspace{.5cm}
 {\it RGZ}^{{\it x2} }_{{\it x4} }=\frac{3}{16}\frac{{\it y4}^{2}\left({\it y2}^{3}+{\it y3}^{3}+{\it y4}^{3}\right)Z^{{\it x2} }}{{\it y2}^{4}}\]}
\end{maplelatex}
 \mapleresult
 \begin{maplelatex}
\mapleinline{inert}{2d}{}{}{\[{\it RGZ}^{{\it x3} }_{{\it x2} }=\frac{3}{16}\frac{\left({\it y2}^{3}+{\it y3}^{3}+{\it y4}^{3}\right)Z^{{\it x3} }}{{\it y2}^{2}}\hspace{1cm}
 {\it RGZ}^{{\it x3} }_{{\it x3} }=-\frac{\left({\it3 y2}^{3}+{\it y3}^{3}+{\it y4}^{3}\right)Z^{{\it x2} }}{{\it 16y2}^{2}}+\frac{{\it 9y4}^{2}Z^{{\it x4} }}{{\it 16 y2} }\]}
\end{maplelatex}
 \mapleresult
\begin{maplelatex}
\mapleinline{inert}{2d}{}{}{\[{\it RGZ}^{{\it x3} }_{{\it x4} }=-\frac{9}{16}\frac{{\it y4}^{2}Z^{{\it x3} }}{{\it y2} }\hspace{2.5cm}
 {\it RGZ}^{{\it x4} }_{{\it x2} }=\frac{3}{16}\frac{\left({\it y2}^{3}+{\it y3}^{3}+{\it y4}^{3}\right)Z^{{\it x4} }}{{\it y2}^{2}}\]}
\end{maplelatex}
 \mapleresult
\begin{maplelatex}
\mapleinline{inert}{2d}{}{}{\[{\it RGZ}^{{\it x4} }_{{\it x3} }=-\frac{9}{16}\frac{{\it y3}^{2}Z^{{\it x4} }}{{\it y2} }\hspace{2.5cm}
 {\it RGZ}^{{\it x4} }_{{\it x4} }=-\frac{3}{16}\frac{\left({\it y2}^{3}+{\it y3}^{3}+{\it y4}^{3}\right)Z^{{\it x2} }}{{\it y2}^{2}}+\frac{9}{16}\frac{{\it y3}^{2}Z^{{\it x3} }}{{\it y2} }\]}
\end{maplelatex}
 \mapleresult
\end{maplegroup}

\bigskip

Putting ${\it RGZ}^{{\it h} }_{{\it i} }=0$, we get a system of algebraic equations. In the case where  $y2^3+y3^3+y4^3=0$, we get the solution $Z^1=t_1$, $Z^2=t_2$ and $Z^3=Z^4=0$ where $t_1,t_2\in \mathbb{R}$. Then,
\begin{equation}\label{nullg}
Z=t_1h_1+t_2h_2.
\end{equation}

\vspace{7pt}

\begin{maplegroup}
\begin{Maple Normal}{
\textbf{ h-curvature  $\overcirc{R}$ of Berwald connection:}}\end{Maple Normal}

\end{maplegroup}
\begin{maplegroup}
\begin{mapleinput}
\mapleinline{active}{1d}{show(RB[i, -h, -j, -k]); }{\[\]}
\end{mapleinput}
\mapleresult
\begin{maplelatex}
\mapleinline{inert}{2d}{}{}{\[
 {\it RB}^{{\it x2} }_{{\it x2 x2 x3} }=\frac{3}{16}\frac{\left({\it y2}^{3}+4{\it y4}^{3}+4{\it y3}^{3}\right){\it y3}^{2}}{{\it y2}^{5}}\hspace{.7cm}
 {\it RB}^{{\it x2} }_{{\it x3 x2 x3} }=-\frac{3}{16}\frac{\left(2{\it y2}^{3}+2{\it y4}^{3}+5{\it y3}^{3}\right){\it y3} }{{\it y2}^{4}}\]}
\end{maplelatex}
 \mapleresult
\begin{maplelatex}
\mapleinline{inert}{2d}{}{}{\[
 {\it RB}^{{\it x2} }_{{\it x4 x2 x3} }=-\frac{9}{16}\frac{{\it y4}^{2}{\it y3}^{2}}{{\it y2}^{4}}\hspace{.7cm}
 {\it RB}^{{\it x3} }_{{\it x2 x2 x3} }=\frac{3}{16}\frac{{\it y2}^{3}-2{\it y3}^{3}-2{\it y4}^{3}}{{\it y2}^{3}}\hspace{.7cm}
 {\it RB}^{{\it x3} }_{{\it x3 x2 x3} }=\frac{9}{16}\frac{{\it y3}^{2}}{{\it y2}^{2}}\]}
\end{maplelatex}
\mapleresult
\begin{maplelatex}
\mapleinline{inert}{2d}{}{}{\[
 {\it RB}^{{\it x3} }_{{\it x4 x2 x3} }=\frac{9}{16}\frac{{\it y4}^{2}}{{\it y2}^{2}}\hspace{.4cm}
 {\it RB}^{{\it x2} }_{{\it x2 x2 x4} }=\frac{3}{16}\frac{\left({\it y2}^{3}+4{\it y4}^{3}+4{\it y3}^{3}\right){\it y4}^{2}}{{\it y2}^{5}}\hspace{.4cm}
 {\it RB}^{{\it x2} }_{{\it x3 x2 x4} }=-\frac{9}{16}\frac{{\it y4}^{2}{\it y3}^{2}}{{\it y2}^{4}}\]}
\end{maplelatex}
 \mapleresult
\begin{maplelatex}
\mapleinline{inert}{2d}{}{}{\[
 {\it RB}^{{\it x2} }_{{\it x4 x2 x4} }=-\frac{3}{16}\frac{\left(2{\it y2}^{3}+5{\it y4}^{3}+2{\it y3}^{3}\right){\it y4} }{{\it y2}^{4}}\hspace{1cm}
 {\it RB}^{{\it x4} }_{{\it x2 x2 x4} }=\frac{3}{16}\frac{{\it y2}^{3}-2{\it y3}^{3}-2{\it y4}^{3}}{{\it y2}^{3}}\]}
\end{maplelatex}
 \mapleresult
\begin{maplelatex}
\mapleinline{inert}{2d}{}{}{\[
 {\it RB}^{{\it x4} }_{{\it x3 x2 x4} }=\frac{9}{16}\frac{{\it y3}^{2}}{{\it y2}^{2}}\hspace{1cm}
 {\it RB}^{{\it x4} }_{{\it x4 x2 x4} }=\frac{9}{16}\frac{{\it y4}^{2}}{{\it y2}^{2}}\hspace{1cm}
 {\it RB}^{{\it x3} }_{{\it x2 x3 x4} }=-\frac{9}{16}\frac{{\it y4}^{2}}{{\it y2}^{2}}\]}
\end{maplelatex}
 \mapleresult
\begin{maplelatex}
\mapleinline{inert}{2d}{}{}{\[
 {\it RB}^{{\it x3} }_{{\it x4 x3 x4} }=\frac{9}{8}\frac{{\it y4} }{{\it y2} }\hspace{1cm}
 {\it RB}^{{\it x4} }_{{\it x2 x3 x4} }=\frac{9}{16}\frac{{\it y3}^{2}}{{\it y2}^{2}}\hspace{1cm}
 {\it RB}^{{\it x4} }_{{\it x3 x3 x4} }=-\frac{9}{8}\frac{{\it y3} }{{\it y2} }\]}
\end{maplelatex}
\end{maplegroup}
\vspace{7pt}

\begin{maplegroup}
\begin{Maple Normal}{
\textbf{ $\overcirc{R}$-nullity vectors}}\end{Maple Normal}

\end{maplegroup}

\begin{maplegroup}
\begin{mapleinput}
\mapleinline{active}{1d}{definetensor(RBW[i, -h, -k] = RB[i, -h, -j, -k]*W[j]); }{\[\]}
\end{mapleinput}
\end{maplegroup}
\begin{maplegroup}
\begin{mapleinput}
\mapleinline{active}{1d}{show(RBW[i, -h, -k]); }{\[\]}
\end{mapleinput}
\mapleresult
\begin{maplelatex}
\mapleinline{inert}{2d}{}{}{\[ {\it RBW}^{{\it x2} }_{{\it x2 x2} }=-\frac{3}{16}\frac{\left({\it y2}^{3}+4{\it y4}^{3}+4{\it y3}^{3}\right){\it y3}^{2}W^{{\it x3} }}{{\it y2}^{5}}-\frac{3}{16}\frac{\left({\it y2}^{3}+4{\it y4}^{3}+4{\it y3}^{3}\right){\it y4}^{2}W^{{\it x4} }}{{\it y2}^{5}}\]}
\end{maplelatex}
\mapleresult
\begin{maplelatex}
\mapleinline{inert}{2d}{}{}{\[ {\it RBW}^{{\it x2} }_{{\it x2 x3} }=\frac{3}{16}\frac{\left({\it y2}^{3}+4{\it y4}^{3}+4{\it y3}^{3}\right)W^{{\it x2} }{\it y3}^{2}}{{\it y2}^{5}}\hspace{.8cm} {\it RBW}^{{\it x2} }_{{\it x2 x4} }=\frac{3}{16}\frac{\left({\it y2}^{3}+4{\it y4}^{3}+4{\it y3}^{3}\right)W^{{\it x2} }{\it y4}^{2}}{{\it y2}^{5}}\]}
\end{maplelatex}
\mapleresult
\begin{maplelatex}
\mapleinline{inert}{2d}{}{}{\[ {\it RBW}^{{\it x2} }_{{\it x3 x2} }=\frac{3}{16}\frac{\left(2{\it y2}^{3}+2{\it y4}^{3}+5{\it y3}^{3}\right){\it y3} W^{{\it x3} }}{{\it y2}^{4}}+\frac{9}{16}\frac{{\it y4}^{2}{\it y3}^{2}W^{{\it x4} }}{{\it y2}^{4}}\hspace{.5cm} {\it RBW}^{{\it x2} }_{{\it x3 x4} }=-\frac{9}{16}\frac{W^{{\it x2} }{\it y4}^{2}{\it y3}^{2}}{{\it y2}^{4}}\]}
\end{maplelatex}
\mapleresult
\begin{maplelatex}
\mapleinline{inert}{2d}{}{}{\[ {\it RBW}^{{\it x2} }_{{\it x3 x3} }=-\frac{3}{16}\frac{\left(2{\it y2}^{3}+2{\it y4}^{3}+5{\it y3}^{3}\right)W^{{\it x2} }{\it y3} }{{\it y2}^{4}}\hspace{1cm} {\it RBW}^{{\it x2} }_{{\it x4 x3} }=-\frac{9}{16}\frac{W^{{\it x2} }{\it y4}^{2}{\it y3}^{2}}{{\it y2}^{4}}\]}
\end{maplelatex}
\mapleresult
\begin{maplelatex}
\mapleinline{inert}{2d}{}{}{\[ {\it RBW}^{{\it x2} }_{{\it x4 x2} }=\frac{9}{16}\frac{{\it y4}^{2}{\it y3}^{2}W^{{\it x3} }}{{\it y2}^{4}}+\frac{3}{16}\frac{\left(2{\it y2}^{3}+5{\it y4}^{3}+2{\it y3}^{3}\right){\it y4} W^{{\it x4} }}{{\it y2}^{4}}\hspace{1cm} {\it RBW}^{{\it x3} }_{{\it x4 x4} }=\frac{9}{8}\frac{W^{{\it x3} }{\it y4} }{{\it y2} }\]}
\end{maplelatex}
\mapleresult
\begin{maplelatex}
\mapleinline{inert}{2d}{}{}{\[ {\it RBW}^{{\it x2} }_{{\it x4 x4} }=-\frac{3}{16}\frac{\left(2{\it y2}^{3}+5{\it y4}^{3}+2{\it y3}^{3}\right)W^{{\it x2} }{\it y4} }{{\it y2}^{4}}\hspace{.5cm}{\it RBW}^{{\it x3} }_{{\it x2 x2} }=-\frac{3}{16}\frac{\left({\it y2}^{3}-2{\it y3}^{3}-2{\it y4}^{3}\right)W^{{\it x3} }}{{\it y2}^{3}}\]}
\end{maplelatex}
\mapleresult
\begin{maplelatex}
\mapleinline{inert}{2d}{}{}{\[ {\it RBW}^{{\it x3} }_{{\it x2 x3} }=\frac{3}{16}\frac{\left({\it y2}^{3}-2{\it y3}^{3}-2{\it y4}^{3}\right)W^{{\it x2} }}{{\it y2}^{3}}+\frac{9}{16}\frac{W^{{\it x4} }{\it y4}^{2}}{{\it y2}^{2}}\hspace{1cm} {\it RBW}^{{\it x3} }_{{\it x2 x4} }=-\frac{9}{16}\frac{W^{{\it x3} }{\it y4}^{2}}{{\it y2}^{2}}\]}
\end{maplelatex}
\mapleresult
\begin{maplelatex}
\mapleinline{inert}{2d}{}{}{\[ {\it RBW}^{{\it x3} }_{{\it x3 x2} }=-\frac{9}{16}\frac{W^{{\it x3} }{\it y3}^{2}}{{\it y2}^{2}}\hspace{.5cm}{\it RBW}^{{\it x3} }_{{\it x3 x3} }=\frac{9}{16}\frac{W^{{\it x2} }{\it y3}^{2}}{{\it y2}^{2}}\hspace{.5cm} {\it RBW}^{{\it x3} }_{{\it x4 x2} }=-\frac{9}{16}\frac{W^{{\it x3} }{\it y4}^{2}}{{\it y2}^{2}}\]}
\end{maplelatex}
\mapleresult
\begin{maplelatex}
\mapleinline{inert}{2d}{}{}{\[ {\it RBW}^{{\it x3} }_{{\it x4 x3} }=\frac{9}{16}\frac{W^{{\it x2} }{\it y4}^{2}}{{\it y2}^{2}}-\frac{9}{8}\frac{{\it y4} W^{{\it x4} }}{{\it y2} }\hspace{1cm}{\it RBW}^{{\it x4} }_{{\it x2 x2} }=-\frac{3}{16}\frac{\left({\it y2}^{3}-2{\it y3}^{3}-2{\it y4}^{3}\right)W^{{\it x4} }}{{\it y2}^{3}}\]}
\end{maplelatex}
\mapleresult
\begin{maplelatex}
\mapleinline{inert}{2d}{}{}{\[ {\it RBW}^{{\it x4} }_{{\it x2 x3} }=-\frac{9}{16}\frac{W^{{\it x4} }{\it y3}^{2}}{{\it y2}^{2}}\hspace{1cm} {\it RBW}^{{\it x4} }_{{\it x2 x4} }=\frac{3}{16}\frac{\left({\it y2}^{3}-2{\it y3}^{3}-2{\it y4}^{3}\right)W^{{\it x2} }}{{\it y2}^{3}}+\frac{9}{16}\frac{W^{{\it x3} }{\it y3}^{2}}{{\it y2}^{2}}\]}
\end{maplelatex}
\mapleresult
\begin{maplelatex}
\mapleinline{inert}{2d}{}{}{\[ {\it RBW}^{{\it x4} }_{{\it x3 x2} }=-\frac{9}{16}\frac{W^{{\it x4} }{\it y3}^{2}}{{\it y2}^{2}}\hspace{.5cm} {\it RBW}^{{\it x4} }_{{\it x3 x3} }=\frac{9}{8}\frac{W^{{\it x4} }{\it y3}} {{\it y2} }\hspace{.5cm} {\it RBW}^{{\it x4} }_{{\it x4 x4} }=\frac{9}{16}\frac{W^{{\it x2} }{\it y4}^{2}}{{\it y2}^{2}}\]}
\end{maplelatex}
\mapleresult
\begin{maplelatex}
\mapleinline{inert}{2d}{}{}{\[ {\it RBW}^{{\it x4} }_{{\it x3 x4} }=\frac{9}{16}\frac{W^{{\it x2} }{\it y3}^{2}}{{\it y2}^{2}}-\frac{9}{8}\frac{{\it y3} W^{{\it x3} }}{{\it y2} }\hspace{1cm}{\it RBW}^{{\it x4} }_{{\it x4 x2} }=-\frac{9}{16}\frac{W^{{\it x4} }{\it y4}^{2}}{{\it y2}^{2}}\]}
\end{maplelatex}
\mapleresult
\end{maplegroup}

\bigskip

Putting ${\it RBW}^{{\it h} }_{{\it ij} }=0$,    we obtain a system of  algebraic equations.  This system has the solution $W^1=t, \, t\in\mathbb{R}$ and  $W^2=W^3=W^4=0$. Then,
\begin{equation}\label{nullb}
W=th_1.
\end{equation}
 Consequently,   (\ref{nullg}) and (\ref{nullb}) lead to $\N_\mathfrak{R} \not\subset \N_{R^\circ}$.

\Section{ Conclusion}

 In this paper, we have mainly achieved  two   objectives:

  $\bullet$ A computational technique  for calculating  the nullity  and kernel vectors, based on the NF-package, has been introduced.

  $\bullet$ Using  this technique,   three  counterexamples have been presented: the first  shows that the two distributions  $\mathrm{Ker}_R$ and  $\N_R$ do not coincide. The second  proves that the nullity distribution $\N_{P^\circ}$  is not completely integrable. The third  shows that the nullity distribution $\N_\mathfrak{R}$  is not a sub-distribution  of  $\N_{R^\circ}$.

%%%%%%%%%%%%%%%%%%%%%%%%%%%%%%%%%%%%%%%%%%%%%%%%%%%%%%%%%%%%%%%%%%%%%%%%%%%%%%%%%%%%%%%%%%%%%%%%%%%%%%%%%%%%%%%%&&&&&&&&&&&&&&&&&&&&&&
\providecommand{\bysame}{\leavevmode\hbox
to3em{\hrulefill}\thinspace}
\providecommand{\MR}{\relax\ifhmode\unskip\space\fi MR }
\providecommand{\MRhref}[2]{%
  \href{http://www.ams.org/mathscinet-getitem?mr=#1}{#2}
} \providecommand{\href}[2]{#2}


\begin{thebibliography}{10}

\bibitem{hbfinsler1}  P. L. Antonelli (Ed.), \emph{Handbook of Finsler
geometry I, II}, Kluwer Acad. publ., 2003.

\bibitem{r101}
P. L. Antonelli, R. Ingarden and M. Matsumoto, \emph{The theory of sprays and Finsler spaces with applications in physics and biology},  Kluwer Acad. Publ., 1993.

\bibitem{Rutz2} P. L. Antonelli,  S. F. Rutz and K. T. Fonseca, \emph{The mathematical theory of endosymbiosis, II: Models of the Fungal Fusion hypothesis},  Nonlinear Anal.,  Real World Appl.,  13 (2012) 2096.


\bibitem{shen2005}
  S. S. Chern and  Z. Shen,  \emph{Riemann-Finsler Geometry}, Singapore: World Scientific, 2005.



\bibitem{r21}
J. Grifone, \emph{Structure presque-tangente et connexions,
\textsc{I}}, Ann.
  Inst. Fourier, Grenoble, \textbf{{22, 1}} (1972), 287--334.
\bibitem{r22}
J. Grifone, \emph{Structure presque-tangente et connexions,
\textsc{II}}, Ann.
  Inst. Fourier, Grenoble, {\textbf{22, 3}} (1972), 291--338.




\bibitem{r27}
J. Klein and A. Voutier, \emph{Formes ext\'{e}rieures
g\'{e}n\'{e}ratrices de
  sprays}, Ann. Inst. Fourier, Grenoble, {\textbf{18, 1}} (1968), 241--260.

  \bibitem{r93}
R. Miron and M. Anastasiei, \emph{The geometry of Lagrange spaces:   Theory and applications}, Kluwer Acad. Publ.,  1994.


\bibitem{gamal} G. G. L. Nashed, \emph{Reissner-nordstr$\ddot{o}$m solutions and energy in teleparallel theory}, Mod. Phys. Lett. A, \textbf{21} (2006), 2241--2250.

\bibitem{Portugal1} R.  Portugal,  S. L.  Sautu,  \emph{Applications  of Maple  to General  Relativity}, Comput. Phys. Commun.,  \textbf{105}  (1997),  233--253.



\bibitem{Rutz3}
S. F. Rutz  and R. Portugal,   \emph{FINSLER: A computer algebra package for Finsler geometries}, Nonlinear Analysis, \textbf{47}  (2001), 6121--6134.

\bibitem{wanas}M. I. Wanas, \emph{On the relation between mass and chage: a pure geometric approach}, Int. J. Geom. Meth. Mod. Phys., \textbf{4} (2007), 373--388.

\bibitem{Nabil.2}
Nabil L. Youssef, \emph{Distribution de nullit\'{e} du tensor de courbure
d'une connexion}, C. R. Acad. Sci. Paris, S\'{e}r. A,
\textbf{290 }(1980), 653--656.

\bibitem{Nabil.1}
Nabil L. Youssef, \emph{Sur les tenseurs de courbure de la connexion de Berwald et ses
distributions de nullit\'{e}.}
Tensor, N. S., \textbf{36} (1982), 275-–280.

\bibitem{ND-Zadeh}
Nabil L. Youssef and S. G. Elgendi, \emph{A note on \lq\lq Sur le noyau de l'op\'{e}rateur de courbure d'une vari\'{e}t\'{e} finsl\'{e}rienne, C. R. Acad. Sci. Paris, s\'er. A, t. 272 (1971), 807-810\rq\rq},  C. R. Math., Ser. I, \textbf{351} (2013), 829--832.
ArXiv: 1305.4498  [math. DG].

\bibitem{CFG}
Nabil~L. Youssef and S. G. Elgendi, \emph{New Finsler   package},  Comput. Phys. Commun., \textbf{185} (2014) 986--997.
    ArXiv: 1306.0875  [math. DG].


\bibitem{ND-cartan}
Nabil~L. Youssef,  A.~Soleiman and S. G. Elgendi, \emph{Nullity distributions associated to
Cartan connection}, Ind. J. Pure Appl. Math., \textbf{45}(2) (2014), 213--238.
ArXiv: 1210.8359  [math. DG].



\end{thebibliography}
\end{document}